\documentclass[12pt]{article}

\usepackage{amsmath,amssymb,amsthm,amscd,a4wide}

\begin{document}

\newcommand{\End}{{\rm{End}\ts}}
\newcommand{\Hom}{{\rm{Hom}}}
\newcommand{\Mat}{{\rm{Mat}}}
\newcommand{\ch}{{\rm{ch}\ts}}
\newcommand{\chara}{{\rm{char}\ts}}
\newcommand{\diag}{ {\rm diag}}
\newcommand{\non}{\nonumber}
\newcommand{\wt}{\widetilde}
\newcommand{\wh}{\widehat}
\newcommand{\ot}{\otimes}
\newcommand{\la}{\lambda}
\newcommand{\La}{\Lambda}
\newcommand{\De}{\Delta}
\newcommand{\al}{\alpha}
\newcommand{\be}{\beta}
\newcommand{\ga}{\gamma}
\newcommand{\Ga}{\Gamma}
\newcommand{\ep}{\epsilon}
\newcommand{\ka}{\kappa}
\newcommand{\vk}{\varkappa}
\newcommand{\si}{\sigma}
\newcommand{\vp}{\varphi}
\newcommand{\de}{\delta}
\newcommand{\ze}{\zeta}
\newcommand{\om}{\omega}
\newcommand{\ee}{\epsilon^{}}
\newcommand{\su}{s^{}}
\newcommand{\hra}{\hookrightarrow}
\newcommand{\ve}{\varepsilon}
\newcommand{\ts}{\,}
\newcommand{\vac}{\mathbf{1}}
\newcommand{\di}{\partial}
\newcommand{\qin}{q^{-1}}
\newcommand{\tss}{\hspace{1pt}}
\newcommand{\Sr}{ {\rm S}}
\newcommand{\U}{ {\rm U}}
\newcommand{\BL}{ {\overline L}}
\newcommand{\BE}{ {\overline E}}
\newcommand{\BP}{ {\overline P}}
\newcommand{\AAb}{\mathbb{A}\tss}
\newcommand{\CC}{\mathbb{C}\tss}
\newcommand{\KK}{\mathbb{K}\tss}
\newcommand{\QQ}{\mathbb{Q}\tss}
\newcommand{\SSb}{\mathbb{S}\tss}
\newcommand{\ZZ}{\mathbb{Z}\tss}
\newcommand{\X}{ {\rm X}}
\newcommand{\Y}{ {\rm Y}}
\newcommand{\Z}{{\rm Z}}
\newcommand{\Ac}{\mathcal{A}}
\newcommand{\Lc}{\mathcal{L}}
\newcommand{\Mc}{\mathcal{M}}
\newcommand{\Pc}{\mathcal{P}}
\newcommand{\Qc}{\mathcal{Q}}
\newcommand{\Tc}{\mathcal{T}}
\newcommand{\Sc}{\mathcal{S}}
\newcommand{\Bc}{\mathcal{B}}
\newcommand{\Ec}{\mathcal{E}}
\newcommand{\Fc}{\mathcal{F}}
\newcommand{\Hc}{\mathcal{H}}
\newcommand{\Uc}{\mathcal{U}}
\newcommand{\Vc}{\mathcal{V}}
\newcommand{\Wc}{\mathcal{W}}
\newcommand{\Ar}{{\rm A}}
\newcommand{\Br}{{\rm B}}
\newcommand{\Ir}{{\rm I}}
\newcommand{\Fr}{{\rm F}}
\newcommand{\Jr}{{\rm J}}
\newcommand{\Or}{{\rm O}}
\newcommand{\GL}{{\rm GL}}
\newcommand{\Spr}{{\rm Sp}}
\newcommand{\Rr}{{\rm R}}
\newcommand{\Zr}{{\rm Z}}
\newcommand{\gl}{\mathfrak{gl}}
\newcommand{\middd}{{\rm mid}}
\newcommand{\Pf}{{\rm Pf}}
\newcommand{\Norm}{{\rm Norm\tss}}
\newcommand{\oa}{\mathfrak{o}}
\newcommand{\spa}{\mathfrak{sp}}
\newcommand{\osp}{\mathfrak{osp}}
\newcommand{\g}{\mathfrak{g}}
\newcommand{\h}{\mathfrak h}
\newcommand{\n}{\mathfrak n}
\newcommand{\z}{\mathfrak{z}}
\newcommand{\Zgot}{\mathfrak{Z}}
\newcommand{\p}{\mathfrak{p}}
\newcommand{\sll}{\mathfrak{sl}}
\newcommand{\agot}{\mathfrak{a}}
\newcommand{\qdet}{ {\rm qdet}\ts}
\newcommand{\Ber}{ {\rm Ber}\ts}
\newcommand{\HC}{ {\mathcal HC}}
\newcommand{\cdet}{ {\rm cdet}}
\newcommand{\tr}{ {\rm tr}}
\newcommand{\str}{ {\rm str}}
\newcommand{\loc}{{\rm loc}}
\newcommand{\Gr}{{\rm G}}
\newcommand{\sgn}{ {\rm sgn}\ts}
\newcommand{\ba}{\bar{a}}
\newcommand{\bb}{\bar{b}}
\newcommand{\bi}{\bar{\imath}}
\newcommand{\bj}{\bar{\jmath}}
\newcommand{\bk}{\bar{k}}
\newcommand{\bl}{\bar{l}}
\newcommand{\hb}{\mathbf{h}}
\newcommand{\Sym}{\mathfrak S}
\newcommand{\fand}{\quad\text{and}\quad}
\newcommand{\Fand}{\qquad\text{and}\qquad}
\newcommand{\OR}{\qquad\text{or}\qquad}

\renewcommand{\theequation}{\arabic{section}.\arabic{equation}}

\newtheorem{thm}{Theorem}[section]
\newtheorem{lem}[thm]{Lemma}
\newtheorem{prop}[thm]{Proposition}
\newtheorem{cor}[thm]{Corollary}
\newtheorem{conj}[thm]{Conjecture}
\newtheorem*{mthma}{Theorem A}
\newtheorem*{mthmb}{Theorem B}

\theoremstyle{definition}
\newtheorem{defin}[thm]{Definition}

\theoremstyle{remark}
\newtheorem{remark}[thm]{Remark}
\newtheorem{example}[thm]{Example}

\newcommand{\bth}{\begin{thm}}
\renewcommand{\eth}{\end{thm}}
\newcommand{\bpr}{\begin{prop}}
\newcommand{\epr}{\end{prop}}
\newcommand{\ble}{\begin{lem}}
\newcommand{\ele}{\end{lem}}
\newcommand{\bco}{\begin{cor}}
\newcommand{\eco}{\end{cor}}
\newcommand{\bde}{\begin{defin}}
\newcommand{\ede}{\end{defin}}
\newcommand{\bex}{\begin{example}}
\newcommand{\eex}{\end{example}}
\newcommand{\bre}{\begin{remark}}
\newcommand{\ere}{\end{remark}}
\newcommand{\bcj}{\begin{conj}}
\newcommand{\ecj}{\end{conj}}

\newcommand{\bal}{\begin{aligned}}
\newcommand{\eal}{\end{aligned}}
\newcommand{\beq}{\begin{equation}}
\newcommand{\eeq}{\end{equation}}
\newcommand{\ben}{\begin{equation*}}
\newcommand{\een}{\end{equation*}}

\newcommand{\bpf}{\begin{proof}}
\newcommand{\epf}{\end{proof}}

\def\beql#1{\begin{equation}\label{#1}}

\title{\Large\bf Characteristic maps for the Brauer algebra}

\author{{A. I. Molev\quad and\quad
N. Rozhkovskaya}}

\date{} 
\maketitle

\vspace{30 mm}

\begin{abstract}
The classical characteristic map associates symmetric functions
to characters of the symmetric groups. There are two natural analogues
of this map involving the Brauer algebra. The first of them relies
on the action of the orthogonal or symplectic group on a space of tensors,
while the second is provided by the action of this group on the symmetric
algebra of the corresponding Lie algebra. We consider
the second characteristic map both in the orthogonal and symplectic case,
and calculate the images of central idempotents
of the Brauer algebra in terms of the Schur polynomials. The calculation
is based on the Okounkov--Olshanski binomial formula for the classical
Lie groups. We also reproduce the hook dimension formulas for representations
of the classical groups by deriving them from the properties of the primitive
idempotents of the symmetric group and the Brauer algebra.
\end{abstract}


\vspace{35 mm}

\noindent
School of Mathematics and Statistics\newline
University of Sydney,
NSW 2006, Australia\newline
alexander.molev@sydney.edu.au

\vspace{7 mm}

\noindent
Department of Mathematics\newline
Kansas State University, USA\newline
rozhkovs@math.ksu.edu

\newpage

\section{Introduction}
\label{sec:int}
\setcounter{equation}{0}

By the classical Schur--Weyl duality,
the natural actions of the
symmetric group $\Sym_m$ and
the general linear group $\GL_N=\GL_N(\CC)$ on the space of tensors
\beql{tenprk}
\underbrace{\CC^{N}\ot\dots\ot\CC^{N}}_m
\eeq
centralize each other. This leads to the multiplicity free decomposition
of the space \eqref{tenprk} as a representation of the group $\Sym_m\times \GL_N$,
\beql{tenprdec}
(\CC^N)^{\ot m}\cong
\bigoplus_{\la\vdash m,\ \ell(\la)\leqslant N}V_{\la}\ot L(\la),
\eeq
where $V_{\la}$ and $L(\la)$
are the respective irreducible representations of $\Sym_m$ and $\GL_N$
associated with a Young diagram $\la$ which contains
$|\la|=m$ boxes, and the number of nonzero rows $\ell(\la)$ does not exceed $N$.

For any element $X\in\End\CC^N$ and $a=1,\dots,m$ we denote by $X_a$ the
corresponding element of the tensor product
\ben
X_a=1^{\ot(a-1)}\ot X\ot 1^{\ot(m-a)}\in \End(\CC^N)^{\ot m}.
\een
An arbitrary element $C$ of the group algebra $\CC[\Sym_m]$
will be regarded as an operator in the space \eqref{tenprk}.
We will identify the symmetric
algebra $\Sr(\gl_N)$ with the algebra of polynomial functions
on the Lie algebra $\gl_N$. If we let $X$ range over $\gl_N$ then the
polynomial function
$X\mapsto \tr\ts C\tss X_1\dots X_m$ with the trace
taken over all $m$ copies of $\End\CC^N$
is a $\GL_N$-invariant element of the
algebra $\Sr(\gl_N)$. This
follows easily by noting that for any matrix $Z\in \GL_N$ we have
\begin{multline}\non
\tr\ts C\ts Z_1 X_1Z_1^{-1}\dots Z_m X_m Z_m^{-1}
=\tr\ts C\ts Z_1\dots Z_m
\tss X_1\dots X_m Z_1^{-1}\dots Z_m^{-1}\\
{}=\tr\ts Z_1^{-1}\dots Z_m^{-1}
\tss C\ts Z_1\dots Z_m
\tss X_1\dots X_m=\tr\ts C\tss X_1\dots X_m,
\end{multline}
where we used the cyclic property of trace and the fact
that the action of the element $C$ commutes with the action
of $\GL_N$.
The algebra of invariants $\Sr(\gl_N)^{\GL_N}$ is isomorphic to
the algebra of symmetric polynomials in $N$ variables. An isomorphism is provided
by the restriction of polynomial functions to
the subspace of diagonal matrices in $\gl_N$.
Hence, the function which takes a diagonal matrix $X$ with eigenvalues
$x_1,\dots,x_N$ to the trace
$\tr\ts C\tss X_1\dots X_m$ is a symmetric polynomial in $x_1,\dots,x_N$.
Thus we can define a linear map
\beql{chara}
\ch:\CC[\Sym_m]\to \CC[x_1,\dots,x_N]^{\Sym_N},\qquad C\mapsto
\frac{1}{m!}\ts\tr\ts C\tss X_1\dots X_m.
\eeq

Given any standard tableau $T$ of shape $\la$ consider the corresponding primitive
idempotent $E_T\in\CC[\Sym_m]$ and calculate its image under the map \eqref{chara}.
To this end, we may assume without loss of generality that the matrix
$X$ is invertible so that $X$ can be viewed as an element of the group $\GL_N$.
The space $E_T(\CC^N)^{\ot m}$ is an irreducible representation of $\GL_N$
isomorphic to $L(\la)$. Therefore, the trace
$\tr\ts E_T X_1\dots X_m$
coincides with the character of the representation $L(\la)$ evaluated at the element $X$.
This value is given by the Weyl character formula so that the trace equals the
Schur polynomial $s_{\la}$ evaluated at the eigenvalues $x_1,\dots,x_N$ of the matrix $X$,
\beql{idemptr}
\tr\ts E_T X_1\dots X_m=s_{\la}(x_1,\dots,x_N).
\eeq
The trace does not depend on the choice of the standard tableau $T$ of shape $\la$ so that
this relation allows us to calculate the image of the irreducible character $\chi_{\la}$
under the map \eqref{chara}. Indeed,
\beql{chila}
\chi_{\la}=\sum_{s\in\Sym_m} \chi_{\la}(s)\cdot s^{-1}
=\frac{m!}{\dim\la}\ts \sum_{T} E_T,
\eeq
where $\dim\la=\dim V_{\la}$ and the second sum is taken over
the standard tableaux $T$ of shape $\la$. Hence
\beql{imchar}
\ch:\chi_{\la}\mapsto s_{\la}(x_1,\dots,x_N).
\eeq
This argument essentially recovers the {\it characteristic map\/}
providing an isomorphism between the algebra generated by the irreducible characters
of the symmetric groups and the algebra of
symmetric functions; cf. \cite[Sec.~I.7]{m:sf}.

Our goal in this paper is to extend the correspondence \eqref{imchar} to a map analogous
to \eqref{chara} involving the Brauer algebra and the respective orthogonal
or symplectic group. Now we suppose that the orthogonal group $\Or_N$ or
symplectic group $\Spr_N$ acts on the space \eqref{tenprk}.
The centralizer of this action
in the endomorphism algebra of the tensor product space
coincides with the homomorphic
image of the Brauer algebra $\Bc_m(\om)$
with the parameter $\om$ specialized to $N$ and $-N$, respectively,
in the orthogonal and symplectic case. This implies the tensor
product decomposition analogous to \eqref{tenprdec},
\beql{dectpro}
(\CC^N)^{\ot m}\cong \bigoplus_{f=0}^{\lfloor m/2\rfloor}
\bigoplus_{\overset{\scriptstyle
\la\vdash m-2f}{\la'_1+\la'_2\leqslant N}}V_{\la}\ot L(\la),
\eeq
where $V_{\la}$ and $L(\la)$
are the respective irreducible representations of $\Bc_m(N)$ and $\Or_N$
associated with the diagram $\la$, and we denote
by $\la'$ the conjugate
diagram so that
$\la'_j$ is the number
of boxes in the column $j$ of $\la$; see \cite{w:cg}.
Given a diagram $\la$ with $\la'_1+\la'_2\leqslant N$ denote by $\la^*$
the diagram obtained from $\la$ by replacing the first column
with the column containing $N-\la'_1$ boxes.
The corresponding representations $L(\la)$ and $L(\la^*)$
of the Lie algebra associated with $\Or_N$ are isomorphic. In what follows we will
only be concerned with the representations $L(\la)$ corresponding to diagrams $\la$
with at most $n$ rows, i.e. $\la'_1\leqslant n$, where $N=2n$ or $N=2n+1$.

Similarly, in the symplectic case with $N=2n$,
\beql{dectprsp}
(\CC^N)^{\ot m}\cong \bigoplus_{f=0}^{\lfloor m/2\rfloor}
\bigoplus_{\overset{\scriptstyle
\la\vdash m-2f}{\la_1\leqslant n}}V_{\la}\ot L(\la'),
\eeq
where $\la_i$ denotes the number of boxes in row $i$ of $\la$,
$V_{\la}$ and $L(\la')$
are the respective irreducible representations of
$\Bc_m(-N)$ and $\Spr_N$
associated with $\la$ and $\la'$; see {\it loc. cit.}

We let $\g_N\subset\gl_N$ denote the orthogonal Lie algebra $\oa_N$ or
symplectic Lie algebra $\spa_N$ which is associated with the
corresponding Lie group $\Gr_N=\Or_N$ or $\Gr_N=\Spr_N$.
We will regard any
element $C$ of the respective
Brauer algebra $\Bc_m(N)$ or $\Bc_m(-N)$ as an operator in the space \eqref{tenprk}.
As with the Lie algebra $\gl_N$, we regard the polynomial function taking
$Y\in\g_N$ to the trace
$\tr\ts C\tss Y_1\dots Y_m$ as an element of the symmetric algebra $\Sr(\g_N)$.
This element is $\Gr_N$-invariant which is verified by the same
calculation as with the corresponding element of $\Sr(\gl_N)$ above.

We will work with a particular presentation of the Lie algebra $\g_N$
so that its Cartan subalgebra consists of diagonal matrices.
Suppose that $y_1,\dots,y_n,-y_1,\dots,-y_n$
are the eigenvalues of a diagonal matrix $Y$
for $N=2n$ and
$y_1,\dots,y_n,-y_1,\dots,-y_n,0$ are the eigenvalues of $Y$ for $N=2n+1$.
Then the function which takes $Y$ to the trace
$\tr\ts C\tss Y_1\dots Y_m$ is a symmetric polynomial in the variables
$y^2_1,\dots,y^2_n$. Thus we get a linear map
\beql{charabcd}
\ch:\Bc_m(\pm N)\to \CC[y^2_1,\dots,y^2_n]^{\Sym_n},\qquad C\mapsto
\frac{1}{m!}\ts\tr\ts C\tss Y_1\dots Y_m.
\eeq

The main result of this paper is the calculation of the image $\ch(\phi_{\la})$
of the normalized central idempotent $\phi_{\la}$ of the Brauer algebra associated
with each partition $\la$ of $m$ satisfying
the respective conditions $\la'_1\leqslant n$ and $\la_1\leqslant n$
in the orthogonal and symplectic case.
We show
that $\ch(\phi_{\la})=0$ if $m$ is odd
and give an explicit formula for the symmetric
polynomial $\ch(\phi_{\la})$ as a linear combination of the Schur polynomials
$s_{\nu}(y^2_1,\dots,y^2_n)$ where $\nu$ runs over partitions of $l$ if $m=2l$.

The starting point of our arguments is the analogue of relation \eqref{idemptr}
for the classical group $\Gr_N$. Namely, suppose that $T$ is a standard tableau of shape
$\la$ and let $Z\in\Gr_N$ be a diagonal matrix such that $\det Z=1$.
We let $E_T$ denote the primitive idempotent of the respective
Brauer algebra $\Bc_m(N)$ or $\Bc_m(-N)$, which we regard as an operator in the
space \eqref{tenprk}. Due to the decompositions
\eqref{dectpro} and \eqref{dectprsp}, the subspace
$E_T(\CC^N)^{\ot m}$ is an irreducible representation of $\Gr_N$
isomorphic to $L(\la)$ in the orthogonal case and to $L(\la')$ in the symplectic case.
Therefore, the trace
$\tr\ts E_T Z_1\dots Z_m$
equals the character of the respective representation $L(\la)$ or $L(\la')$ so that
\beql{idemptrbd}
\tr\ts E_T Z_1\dots Z_m=\chi^{\oa_N}_{\la}(z_1,\dots,z_n)
\eeq
in the orthogonal case, and
\beql{idemptrc}
\tr\ts E_T Z_1\dots Z_m=\chi^{\spa_N}_{\la'}(z_1,\dots,z_n)
\eeq
in the symplectic case, where we denote by
$z^{}_1,\dots,z^{}_n,z^{-1}_1,\dots,z^{-1}_n$
the eigenvalues of $Z$
for $N=2n$ and by $z^{}_1,\dots,z^{}_n,z^{-1}_1,\dots,z^{-1}_n, 1$
the eigenvalues of $Z$ for $N=2n+1$. Explicit expressions for the characters
are well known; they are implied by the Weyl character formula
and can be found e.g. in \cite{oo:ss2}.
Although relations
\eqref{idemptrbd} and \eqref{idemptrc} are analogous to \eqref{idemptr},
note a principal difference with the case of $\GL_N$. In that case,
the matrix $X$ in \eqref{idemptr}
could be treated both as an element of the group $\GL_N$ and
as an element of the Lie algebra $\gl_N$. In contrast, the passage from the group $\Gr_N$
to the Lie algebra $\g_N$ requires an additional step.
To derive explicit formulas
for the images of central idempotents of the Brauer algebra under the map
\eqref{charabcd} we use the Okounkov--Olshanski binomial
formula \cite[Theorem~1.2]{oo:ss2}. This allows us
to express the characters occurring in
\eqref{idemptrbd} and \eqref{idemptrc} in terms of the Schur polynomials
in $y_1^2,\dots,y_n^2$.

More precisely, by analogy with \eqref{chila} set
\beql{philao}
\phi_{\la}=\frac{1}{D(\la)}\ts \sum_{T} E_T
\eeq
in the orthogonal case, and
\beql{philasp}
\phi_{\la}=\frac{1}{D(\la')}\ts \sum_{T} E_T
\eeq
in the symplectic case, where $D(\la)=\dim L(\la)$ and the sums are taken
over standard tableaux $T$ of shape $\la$.
Our main result (see Theorem~\ref{thm:mainch} below) states
that $\ch(\phi_{\la})=0$ unless $m$ is even, $m=2l$.
In this case,
\beql{mainf}
\ch(\phi_{\la})
=\sum_{\nu\vdash l}\frac{s_{\nu}(y_1^2,\dots,y_n^2)}{C(\nu)}
\sum_{\mu\subseteq\la}(-1)^{|\mu|}\ts
\frac{s_{\nu}(a_{\rho}|\ts a)}{H(\mu)\ts H(\la/\mu)},
\eeq
where $\rho=\mu$ and $\rho=\mu'$ in the orthogonal and symplectic case, respectively,
and we use the following notation.
For any
skew diagram $\theta$ we denote by
$\dim\theta$
the number of standard $\theta$-tableaux
with
entries in $\{1,2,\dots,|\theta|\}$ and set
\beql{prhook}
H(\theta)=\frac{|\theta|!}{\dim\theta}.
\eeq
If $\theta$ is normal (nonskew), then
$H(\theta)$ coincides with the product of
the hooks of $\theta$ due to the hook formula.
Furthermore, the constant $C(\nu)=C_{\g_N}(\nu)$
equals the inner square of the Schur polynomial
$s_{\nu}(y_1^2,\dots,y_n^2)$ with respect to
an invariant inner product on $\Sr(\g_N)$ (see \cite[Proposition~5.3]{oo:ss2})
and it is given by
\beql{cnu}
C(\nu)=\prod_{(i,j)\in \nu} 2\tss\big(n+j-i\big)\big(N-1+2\tss(j-i+\ve)\big),
\eeq
where
\ben
\ve=\begin{cases}0\qquad&\text{for}\quad \g_N=\oa_{2n},\\
1/2\qquad&\text{for}\quad \g_N=\oa_{2n+1},\\
1\qquad&\text{for}\quad \g_N=\spa_{2n}.
\end{cases}
\een
Finally, by $s_{\nu}(x\ts |\ts a)$ we denote the {\it double\/} (or
{\it factorial\/}) {\it Schur polynomial\/}
in the variables $x=(x_1,\dots,x_n)$ associated with the particular parameter sequence
$a=(a_i\ts |\ts i\in\ZZ)$ with $a_i=(\ve+i-1)^2$. The polynomial $s_{\nu}(x\ts|\ts a)$
is symmetric in $x_1,\dots,x_n$ and it
can be given by several equivalent formulas; see e.g. \cite[Sec.~I.3]{m:sf}
(note that the
sequence $a$ there corresponds to our sequence $-a$).
In particular,
\beql{doubsf}
s_{\nu}(x\ts|\ts a)=
\sum_{T}\prod_{\al\in\nu} (x_{T(\al)}-a_{T(\al)+c(\al)}),
\eeq
summed over semistandard $\nu$-tableaux $T$ with entries in
$\{1,\dots,n\}$, where $c(\al)=j-i$ denotes the content of the box $\al=(i,j)$.
For any partition $\mu$ with at most $n$ parts we denote by $a_{\mu}$
the $n$-tuple
\beql{amun}
a_{\mu}=(a_{\mu_1+n},a_{\mu_2+n-1},\dots, a_{\mu_n+1}).
\eeq

We demonstrate below (Sec.~\ref{subsec:sa}) that the second sum in \eqref{mainf}
simplifies if $\la$ is a row or column diagram. However, we do not know whether
shorter expressions for this sum exist for arbitrary $\la$.

The relation \eqref{mainf} is remarkably similar
to the main result of Nazarov's paper \cite[Theorem~3.4]{n:ce}
involving a different kind of characteristic maps.

As a consequence of our approach, we also demonstrate that the well-known hook
dimension formulas for representations of the classical groups can be obtained
directly from the properties of the primitive
idempotents of the symmetric group and the Brauer algebra
via relations \eqref{idemptr}, \eqref{idemptrbd} and
\eqref{idemptrc} implied by the Schur--Weyl duality.

We believe it is
possible to extend the main results of the paper to the case where $\la$
is a partition of $m-2f$ with $f\geqslant 1$ and to calculate
the images $\ch(\phi_{\la})$ of the associated central idempotents $\phi_{\la}$
under the map \eqref{charabcd}.
This should involve more
complicated combinatorics of paths in the Bratteli diagram corresponding
to the Brauer algebras.

\medskip
We are grateful to Arun Ram for valuable discussions and for providing us
with the preprints \cite{drv:ad} and \cite{drv:ada}.
The second author would like to thank the School of Mathematics and Statistics
of the University of Sydney for warm hospitality during her visit.

\section{Idempotents in the Brauer algebra}
\label{sec:bra}
\setcounter{equation}{0}

Let $m$ be a positive integer and $\om$ an indeterminate.
An $m$-diagram $d$ is a collection
of $2m$ dots arranged into two rows with $m$ dots in each row
connected by $m$ edges such that any dot belongs to only one edge.
The product of two $m$-diagrams $d_1$ and $d_2$ is determined by
placing $d_1$ above $d_2$ and identifying the vertices
of the bottom row of $d_1$ with the corresponding
vertices in the top row of $d_2$. Let $s$ be the number of
closed loops obtained in this placement. The product $d_1d_2$ is given by
$\om^{\tss s}$ times the resulting diagram without loops.
The {\it Brauer algebra\/} $\Bc_m(\om)$ \cite{b:aw} is defined as the
$\CC(\om)$-linear span of the $m$-diagrams with the
multiplication defined above.
The dimension of the algebra is $1\cdot 3\cdots (2m-1)$.

For $1\leqslant a<b\leqslant m$ denote by $s_{a\tss b}$
and $\ep_{a\tss b}$ the respective diagrams
of the form

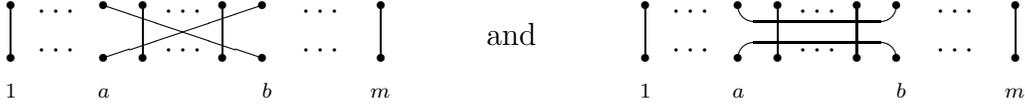
\begin{center}
\begin{picture}(400,60)
\thinlines

\put(10,20){\circle*{3}}
\put(45,20){\circle*{3}}
\put(60,20){\circle*{3}}
\put(90,20){\circle*{3}}
\put(105,20){\circle*{3}}
\put(150,20){\circle*{3}}

\put(10,40){\circle*{3}}
\put(45,40){\circle*{3}}
\put(60,40){\circle*{3}}
\put(90,40){\circle*{3}}
\put(105,40){\circle*{3}}
\put(150,40){\circle*{3}}

\put(10,20){\line(0,1){20}}
\put(60,20){\line(0,1){20}}
\put(45,20){\line(3,1){60}}
\put(45,40){\line(3,-1){60}}
\put(90,20){\line(0,1){20}}
\put(150,20){\line(0,1){20}}

\put(20,20){$\cdots$}
\put(20,35){$\cdots$}
\put(68,20){$\cdots$}
\put(68,35){$\cdots$}
\put(120,20){$\cdots$}
\put(120,35){$\cdots$}

\put(8,5){\scriptsize $1$ }
\put(43,5){\scriptsize $a$ }
\put(105,5){\scriptsize $b$ }
\put(146,5){\scriptsize $m$ }

\put(190,25){\text{and}}

\put(250,20){\circle*{3}}
\put(285,20){\circle*{3}}
\put(300,20){\circle*{3}}
\put(330,20){\circle*{3}}
\put(345,20){\circle*{3}}
\put(390,20){\circle*{3}}

\put(250,40){\circle*{3}}
\put(285,40){\circle*{3}}
\put(300,40){\circle*{3}}
\put(330,40){\circle*{3}}
\put(345,40){\circle*{3}}
\put(390,40){\circle*{3}}

\put(250,20){\line(0,1){20}}
\put(300,20){\line(0,1){20}}
\put(315,20){\oval(60,12)[t]}
\put(315,40){\oval(60,12)[b]}
\put(330,20){\line(0,1){20}}
\put(390,20){\line(0,1){20}}

\put(260,20){$\cdots$}
\put(260,35){$\cdots$}
\put(308,20){$\cdots$}
\put(308,35){$\cdots$}
\put(360,20){$\cdots$}
\put(360,35){$\cdots$}

\put(248,5){\scriptsize $1$ }
\put(283,5){\scriptsize $a$ }
\put(345,5){\scriptsize $b$ }
\put(386,5){\scriptsize $m$ }

\end{picture}
\end{center}

\noindent
and set $s_a=s_{a\ts a+1}$ and $\ee_a=\ee_{a\ts a+1}$ for $a=1,\dots,m-1$.
The subalgebra of $\Bc_m(\om)$ generated over $\CC$
by the elements $s_{a\ts b}$
is isomorphic to the group algebra of the symmetric group
$\CC[\Sym_m]$ so that $s_{a\tss b}$
will be identified with the transposition $(a\ts b)$.
The Brauer algebra $\Bc_{m-1}(\om)$ will be regarded
as a natural subalgebra of $\Bc_m(\om)$.

The {\it Jucys--Murphy elements\/} $x_1,\dots,x_m$
for the Brauer algebra $\Bc_m(\om)$ are given by the formulas
\beql{jmdef}
x_b=\frac{\om-1}{2}+\sum_{a=1}^{b-1}(s_{a\tss b}-\ep_{a\tss b}),
\qquad b=1,\dots,m;
\eeq
see \cite{lr:rh} and \cite{n:yo}, where, in particular,
the eigenvalues
of the $x_b$ in irreducible representations were calculated.
We have followed \cite{n:yo} to include the shift by $(\om-1)/2$ in the definition
to simplify the formulas below.
The element $x_m$ commutes
with the subalgebra $\Bc_{m-1}(\om)$. This implies that
the elements $x_1,\dots,x_m$
of $\Bc_m(\om)$ pairwise commute. A complete set of pairwise orthogonal primitive
idempotents for the Brauer algebra can be constructed
with the use of these elements.
Suppose that $\la$ is a partition of $m$.
We will identify partitions with their diagrams
so that if the parts of $\la$ are $\la_1,\la_2,\dots$ then
the corresponding diagram is
a left-justified array of rows of unit boxes containing
$\la_1$ boxes in the top row, $\la_2$ boxes
in the second row, etc.
The box in row $i$ and column $j$ of a diagram
will be denoted as the pair $(i,j)$.
A {\it standard $\la$-tableau\/} is a sequence
$T=(\La_1,\dots,\La_m)$ of diagrams such that
for each $r=1,\dots,m$ the diagram $\La_r$ is obtained
from $\La_{r-1}$ by adding one box,
where we set $\La_0=\varnothing$ (the empty diagram) and $\La_m=\la$.
Equivalently, $T$ will be viewed as the array obtained by writing
$r\in\{1,\dots,m\}$ into the box of the diagram $\la$
which is added to the diagram $\La_{r-1}$
to get $\La_r$.
To each standard tableau $T$ we associate the corresponding
sequence of {\it contents\/} $(c_1,\dots,c_m)$, $c_a=c_a(T)$,
where
\beql{conbr}
c_a=\frac{\om-1}{2}+j-i
\eeq
if $\La_a$ is obtained by adding the box $(i,j)$ to $\La_{a-1}$.
The primitive idempotents $E_{T}=E^{\lambda}_{T}$
can now be defined by
the following recurrence formula
(we omit the superscripts
indicating the diagrams since they are determined by
the standard tableaux).
Set $\mu=\La_{m-1}$ and consider the standard $\mu$-tableau
$U=(\La_1,\dots,\La_{m-1})$ so that $U$ can be viewed as the
tableau obtained from $T$ by
removing the box containing $m$. Then
\beql{jmform}
E_{T}=E_{U}\ts \frac{u-c_m}
{u-x_m}\ts
\Big|^{}_{u=c_m}.
\eeq
Since the Brauer algebra is finite-dimensional,
the fraction involving $x_m$ on the right hand side reduces
to a polynomial in $x_m$
whose coefficients are rational functions in $u$.
They turn out to be well-defined at $u=c_m$ and relation
\eqref{jmform} states that the evaluation of the right hand side yields $E_{T}$.
This relation is essentially a version of the well-known Jucys--Murphy formula;
see \cite{im:fp} and \cite{m:fp} for more details.

\section{Hook dimension formulas for classical groups}
\label{sec:hf}
\setcounter{equation}{0}

Consider the vector space $\CC^N$
with its canonical basis $e_1,\dots,e_N$.
We will be using the involution on the set
of indices $\{1,\dots,N\}$ defined by $i\mapsto i^{\tss\prime}=N-i+1$
and equip the space $\CC^N$ with the following nondegenerate symmetric
or skew-symmetric bilinear form
\beql{formg}
\langle e_i,e_j\rangle=g_{ij},
\eeq
where $N=2n$ is even in the skew-symmetric case, and
\beql{formgss}
g_{ij}=\begin{cases}
\phantom{\ve_i\ts}\de_{i\ts j^{\tss\prime}}\quad&
\text{in the symmetric case,}\\
\ve_i\ts\de_{i\ts j^{\tss\prime}}\quad&\text{in the skew-symmetric case,}
\end{cases}
\eeq
with $\ve_i=1$ for $i=1,\dots,n$ and
$\ve_i=-1$ for $i=n+1,\dots,2n$.

The classical
group $\Gr_N=\Or_N$ or $\Gr_N=\Spr_{N}$ is defined as the group
of complex matrices
preserving the respective symmetric or skew-symmetric form \eqref{formg},
\ben
\Gr_N=\{Z\in \Mat_N(\CC)\ |\ Z^t G\ts Z=G\},\qquad G=[g_{ij}].
\een

Observe that if $Z=1$ is the identity matrix, then the values provided by
the expressions \eqref{idemptrbd} and \eqref{idemptrc} coincide
with the dimensions $\dim L(\la)$ and $\dim L(\la')$ of the respective
representations of the groups $\Or_N$ and $\Spr_N$. It is well-known by
El~Samra and King~\cite{ek:di} that these dimensions are given by
the {\it hook formulas}. We will consider the images of the idempotents $E_T$
under the action of the Brauer algebra in the tensor space \eqref{tenprk}
and calculate their partial traces with respect to the $m$-th copy of
the vector space $\CC^N$.
In particular, this will provide another proof of the hook dimension formulas of \cite{ek:di};
see also \cite{w:qg}. To make our arguments clearer, we
first go over a technically simpler case of $\GL_N$
to reproduce Robinson's formula; see e.g.~\cite{ek:di}
and \cite[Sec.~I.3, Example~4]{m:sf}
for other proofs.

\subsection{Dimension formulas for $\GL_N$}
\label{subsec:hd}

The symmetric group version of the recurrence relation \eqref{jmform}
takes exactly the same form \cite{m:fp} with the respective definitions
of the objects associated with $\Sym_m$.
Here, as above, $T$ is a standard tableau of shape $\la\vdash m$ and $U$
is obtained from $T$ by removing the box occupied by $m$. The content $c_a=c_a(T)$
of the box $(i,j)$ of $T$ occupied by $a$ is now found by $c_a=j-i$ and
the Jucys--Murphy elements $x_a$ are now given by $x_a=s_{1\ts a}+\dots+s_{a-1\ts a}$;
cf. \eqref{jmdef}.
In the group algebra we have the relation $s_{m-1}x_m=x_{m-1}s_{m-1}+1$
which implies
\beql{smmone}
s_{m-1}+\frac{1}{u-x_m}=(u-x_{m-1})\ts s_{m-1}\ts \frac{1}{u-x_m}.
\eeq
Hence,
\ben
\bal
\frac{1}{u-x_m}&=s_{m-1}\ts \frac{1}{u-x_{m-1}}\ts\Big(s_{m-1}+\frac{1}{u-x_m}\Big)\\
{}&=s_{m-1}\ts \frac{1}{u-x_{m-1}}\ts s_{m-1}+
s_{m-1}\ts\frac{1}{(u-x_{m-1})(u-x_m)}.
\eal
\een
Therefore, applying \eqref{smmone} once again we come to the identity
\beql{frixm}
\frac{1}{u-x_m}=s_{m-1}\ts \frac{1}{u-x_{m-1}}\ts s_{m-1}+
\frac{1}{u-x_{m-1}}\ts\Big(s_{m-1}+\frac{1}{u-x_m}\Big)\ts\frac{1}{u-x_{m-1}}.
\eeq

Consider the action of the symmetric group
$\Sym_m$ in the space \eqref{tenprk} so that the image of the element $s_{ab}\in \Sym_m$
with $a<b$ is found by $s_{ab}\mapsto P_{ab}$,
\beql{pdef}
P_{ab}=\sum_{i,j=1}^N 1^{\ot(a-1)}\ot e_{ij}
\ot 1^{\ot(b-a-1)}\ot e_{ji}\ot 1^{\ot(m-b)},
\eeq
where the $e_{ij}\in\End\CC^{N}$ denote the standard matrix units.
From now on we use this action
to regard elements of the group algebra $\CC[\Sym_m]$
as elements of the algebra $\End\big((\CC^{N})^{\ot m}\big)$ which is
naturally identified
with the tensor product of the endomorphism algebras,
\beql{tenpendrk}
\End\big((\CC^{N})^{\ot m}\big)\cong\underbrace{\End\CC^{N}\ot\dots\ot\End\CC^{N}}_m.
\eeq

The trace map $\tr:\End\CC^{N}\to\CC$ is defined in a usual way as a linear map
taking the matrix unit $e_{ij}$ to $\de_{ij}$.
For each $a=1,\dots,m$ we will consider the partial trace
$\tr_a$ as a linear map $(\End\CC^{N})^{\ot m}\to (\End\CC^{N})^{\ot (m-1)}$
applied to the $a$-th copy of the endomorphism algebra. Note that
$\tr_a(P_{ab})=1$.
Furthermore, since
\ben
\tr_m \Big(s_{m-1}\ts \frac{1}{u-x_{m-1}}\ts s_{m-1}\Big)=\tr^{}_{m-1}
\Big(\frac{1}{u-x_{m-1}}\Big),
\een
we get a recurrence relation for the rational functions
\beql{amuf}
A_m(u)=\tr_m\Big(\frac{1}{u-x_m}\Big)
\eeq
in the form
\ben
A_m(u)=A_{m-1}(u)+\frac{1}{(u-x_{m-1})^2}\ts \big(1+A_m(u)\big),
\een
that is,
\ben
A_m(u)=\frac{(u-x_{m-1})^2}{(u-x_{m-1})^2-1}\ts A_{m-1}(u)+\frac{1}{(u-x_{m-1})^2-1}.
\een
Since for $m=1$ we have $A_1(u)=N/u$, solving the recurrence relation we find that
\ben
A_m(u)=\frac{u+N}{u}\prod_{a=1}^{m-1}\ts\frac{(u-x_a)^2}{(u-x_a)^2-1}-1.
\een

This calculation and the recurrence formula \eqref{jmform} allow us to
find the partial trace $\tr_m E_T$ of the idempotent $E_T$ regarded as
an element of the algebra \eqref{tenpendrk}. By the properties
of the Jucys--Murphy elements,
\beql{xiet}
x_a\ts E_{U}=E_{U}\ts x_a=c_a\ts E_{U}, \qquad a=1,\dots,m-1.
\eeq
Therefore,
\ben
\tr_m E_T=E_{U}\ts \Big[(u-c_m)\ts A_m(u)\Big]^{}_{u=c_m}
=(N+c_m)\ts E_{U}\ts \Big[\ts\frac{u-c_m}{u}\prod_{a=1}^{m-1}\ts
\frac{(u-c_a)^2}{(u-c_a)^2-1}\ts\Big]_{u=c_m}.
\een
The evaluation of the rational function in $u$ is well-defined and it depends
only on the shape $\mu$ of the standard tableau $U$ but does not depend on $U$.
The result of the evaluation is easily calculated (cf. \cite{m:fp}); it gives
\beql{treth}
\tr_m E_T=(N+c_m)\ts \frac{H(\mu)}{H(\la)}\ts E_{U},
\eeq
where $H(\mu)$ and $H(\la)$ are the products of hooks of the diagrams $\mu$ and $\la$.
By \eqref{idemptr}, the dimension of the
irreducible representation $L(\la)$ of $\GL_N$ is found as the trace of $E_T$
taken over all $m$ copies of the endomorphism space $\End\CC^N$. Hence, applying
\eqref{treth}
we arrive at the well-known Robinson formula for this dimension.
If $\la$ is a partition with at most $N$ parts,
the dimension of the irreducible representation $L(\la)$ of $\GL_N$ is given by
\ben
\dim L(\la)=\frac{1}{H(\la)}\prod_{(i,j)\in\la}(N+j-i).
\een

\subsection{Dimension formulas for $\Or_N$ and $\Spr_N$}
\label{subsec:os}

To prove analogues of the hook dimension formula for
the orthogonal and symplectic groups, consider the recurrence
relation \eqref{jmform} in the Brauer algebra $\Bc_m(\om)$.
The starting point will be an analogue of the identity \eqref{frixm} for $\Bc_m(\om)$
given in the next lemma. This identity goes back to \cite[Sec.~4.1]{n:yo} where
it is proved for the degenerate affine Wenzl algebra and used in the description
of the center of that algebra; see also \cite{drv:ad} for generalizations to
the affine BMW algebras.
Recall that now the Jucys--Murphy elements $x_a$ and the contents $c_a$ are
defined by \eqref{jmdef} and \eqref{conbr}.

\ble\label{lem:idium}
We have the identity of rational functions in $u$,
\ben
\bal
\frac{1}{u-x_m}&=
s_{m-1}\ts\frac{1}{u-x_{m-1}}\ts s_{m-1}
\ts +\ts
\frac{1}{u-x_{m-1}}\ts s_{m-1}\ts \frac{1}{u-x_{m-1}}
 \\
& \quad  +\frac{1}{(u-x_{m-1})^2}\ts\frac{1}{u-x_{m}}
\ts -\ts
 \frac{1}{u-x_{m-1}}\ts \ep_{m-1}\ts \frac{1}{(u+x_{m-1})(u-x_{m-1})}
 \ts\\
&\quad -\frac{1}{u+x_{m-1}}\ts\ep_{m-1}\ts \frac{1}{u+x_{m-1}}
\ts+\ts
\frac{1}{(u+x_{m-1})(u-x_{m-1})}\ts\ep_{m-1}\ts \frac{1}{u+x_{m-1}}
\\
&\quad -\frac{1}{u+x_{m-1}}\ts\ep_{m-1}\ts \frac{1}{u-x_{m-1}}
\ts\ep_{m-1}\ts \frac{1}{u+x_{m-1}}.
\eal
\een
\ele

\bpf
Note the following relations in $\Bc_m(\om)$ satisfied by the
Jucys--Murphy elements:
\begin{align}\label{ex}
\ep_{m-1}\,x_m&=-\ep_{m-1} x_{m-1},\\
\label{sx}
s_{m-1}x_m&=x_{m-1}s_{m-1}+1-\ep_{m-1}.
\end{align}
By \eqref{sx} we have $s_{m-1}(u-x_m)=(u-x_{m-1})\ts s_{m-1} -(1-\ep_{m-1})$.
Multiply both sides of this relation by $(u-x_{m-1})^{-1}$ from the left and by
$(u-x_{m})^{-1}$ from the right and rearrange to get
\beql{smoux}
s_{m-1}\ts\frac{1}{u-x_{m}}=\frac{1}{u-x_{m-1}}\ts s_{m-1}
+\frac{1}{u-x_{m-1}}\ts (1-\ep_{m-1})\ts \frac{1}{u-x_m}
\eeq
which implies
\ben
\bal
\frac{1}{u-x_m}=s_{m-1} \ts\frac{1}{u-x_{m-1}}\ts s_{m-1}&+s_{m-1}\ts
\frac{1}{(u-x_{m-1})(u-x_m)}\\
{}&-s_{m-1}\ts
\frac{1}{u-x_{m-1}}\ts\ep_{m-1}\ts\frac{1}{u-x_m}.
\eal
\een
The desired identity will follow after rewriting the second and third terms
on the right hand side with the use of
\eqref{ex}, \eqref{smoux} and the property that
the elements $x_{m-1}$ and $x_m$ commute. The second term takes the form
\begin{align*}
s_{m-1}\ts\frac{1}{(u-x_{m})(u-x_{m-1})}&
=\frac{1}{u-x_{m-1}}\ts s_{m-1}\ts\frac{1}{u-x_{m-1}}
+ \frac{1}{(u-x_{m-1})^2(u-x_m)}\\
{}&- \frac{1}{u-x_{m-1}}\ts\ep_{m-1}\ts\frac{1}{(u+x_{m-1})(u-x_{m-1})},
\end{align*}
while for the third term we have
\begin{align*}
&s_{m-1}\ts\frac{1}{u-x_{m-1}}\ts\ep_{m-1}\ts\frac{1}{u-x_m} =\\
&=\Big(\ts\frac{1}{u-x_{m}}\ts s_{m-1} - \frac{1}{(u-x_{m})(u-x_{m-1})}
+\frac{1}{u-x_{m}}\ts\ep_{m-1}\ts\frac{1}{u-x_{m-1}}
\Big)\ts\ep_{m-1}\ts\frac{1}{u-x_m}\\
&=\frac{1}{u+x_{m-1}}\ts\ep_{m-1}\ts\frac{1}{u+x_{m-1}}
 - \frac{1}{(u-x_{m-1})(u+x_{m-1})}\ts\ep_{m-1}\ts\frac{1}{u+x_{m-1}}\\
 &\quad +\frac{1}{u+x_{m-1}}\ts\ep_{m-1}\ts
 \frac{1}{u-x_{m-1}}\ts\ep_{m-1}\ts\frac{1}{u+x_{m-1}},
\end{align*}
completing the proof.
\epf

Now consider the natural action of the orthogonal group $\Or_N$ in the space
of tensors \eqref{tenprk} and the commuting action of the Brauer algebra $\Bc_m(N)$
so that the parameter $\om$ is specialized to $N$.
The action of $\Bc_m(N)$ in the space \eqref{tenprk}
is defined by the assignments
\beql{braact}
s_{ab}\mapsto P_{ab},\qquad \ep_{ab}\mapsto Q_{ab},\qquad
a< b,
\eeq
where $P_{ab}$ is defined in \eqref{pdef}, and
\beql{qdefo}
Q_{ab}=\sum_{i,j=1}^N 1^{\ot(a-1)}\ot e_{ij}
\ot 1^{\ot(b-a-1)}\ot e_{i'j'}\ot 1^{\ot(m-b)}.
\eeq
Note that $\tr_a(Q_{ab})=1$ for $1\leqslant a<b\leqslant m$, and
for any element $X\in\End\CC^N$
we have the property
$Q_{ab}\ts X_a Q_{ab}=\tr(X)\ts Q_{ab}$.
Now we use Lemma~\ref{lem:idium} and regard elements of the algebra $\Bc_m(N)$
as elements of the algebra \eqref{tenpendrk}. Define the functions $A_m(u)$ by
the same formula \eqref{amuf} as for the symmetric group,
but with the new definition \eqref{jmdef} of the Jucys--Murphy elements.
Calculating the partial trace $\tr_m$ on both sides of
the identity of Lemma~\ref{lem:idium}
we get the recurrence relation
\ben
\bal
A_m(u)&=A_{m-1}(u) +\frac{1}{(u-x_{m-1})^2}+\frac{1}{(u-x_{m-1})^2}\ts A_m(u)-
\frac{1}{(u-x_{m-1})^2(u+x_{m-1})}\\
{}&-\frac{1}{(u+x_{m-1})^2} +\frac{1}{(u+x_{m-1})^2(u-x_{m-1})} -
\frac{1}{(u+x_{m-1})^2}\ts A_{m-1}(u)
\eal
\een
which simplifies to
\ben
\bal
\Big(1- \frac{1}{(u-x_{m-1})^2}\Big)A_m(u)
 =
\Big(1- \frac{1}{(u+x_{m-1})^2}\Big)A_{m-1}(u)
+\frac{2\ts(2\tss u-1)\ts x_{m-1}}{(u-x_{m-1})^2(u+x_{m-1})^2}.
\eal
\een
For $m=1$ we have $A_1(u)=N/(u-c_1)$, where $c_1=(N-1)/2$
and the relation is easily solved by
using the substitution
\ben
A_m(u)=
\wt{A}_m(u)- \frac{2u-1}{2u}.
\een
The solution reads (cf. closely related calculations in \cite{drv:ada}):
\ben
A_m(u)=\prod_{a=1}^{m-1}\frac{(u-x_a)^2}{(u-x_a)^2-1}
\prod_{a=1}^{m-1}\frac{(u+x_a)^2-1}{(u+x_a)^2}
\Big(\frac{2\tss N}{2\tss u-N+1}+\frac{2\tss u-1}{2\tss u}\Big)-\frac{2\tss u-1}{2\tss u}.
\een

For any diagram $\la$ with $\la'_1\leqslant n$ set
\beql{dimf}
D(\la)=\frac{1}{H(\la)}\prod_{(i,j)\in\la}\big(N-1+d(i,j)\big),
\eeq
where $H(\la)$ is the product of hooks of $\la$ (see \eqref{prhook}) and
\ben
d(i,j)=\begin{cases} \la_i+\la_j-i-j+1&\qquad
\text{if} \quad i\leqslant j,\\
-\la'_i-\la'_j+i+j-1&\qquad\text{if} \quad i> j.
\end{cases}
\een
Let $T$ be a standard tableau of shape $\la$. Denote by $U$ the standard tableau
obtained from $T$ by deleting the box occupied by $m$ and let $\mu$ be the shape of $U$.
As with the group algebra $\CC[\Sym_m]$ in Sec.~\ref{subsec:hd},
the recurrence formula \eqref{jmform} allows us to
find the partial trace $\tr_m E_T$ of the idempotent $E_T$ regarded as
an element of the algebra \eqref{tenpendrk}.
The following proposition also recovers the hook dimension formula
\cite[(3.28)]{ek:di}. Note that it is given there in an equivalent form which amounts
to a change in the definition of $d(i,j)$: the inequalities $i\leqslant j$ and $i>j$
are respectively replaced by $i\geqslant j$ and $i<j$.

\bpr\label{prop:parto}
We have the relation
\beql{trete}
\tr_m E_T=E_U\ts \frac{D(\la)}{D(\mu)}.
\eeq
Moreover,
the dimension of the irreducible representation $L(\la)$ of $\Or_N$
equals $D(\la)$.
\epr

\bpf
We have
\ben
\tr_m E_T=E_U\Big[(u-c_m)A_m(u)\Big]_{u=c_m}
\een
and using the above formula for $A_m(u)$ we get
\ben
\tr_m E_T=E_U\Big[\frac{u-c_m}{u-c_1}
\prod_{a=1}^{m-1}\frac{(u-c_a)^2}{(u-c_a)^2-1}
\Big]_{u=c_m}\prod_{a=1}^{m-1}\frac{(c_m+c_a)^2-1}{(c_m+c_a)^2}
\Big(N+\frac{(2\tss c_m-1)(c_m-c_1)}{2\tss c_m}\Big).
\een
As we found in Sec.~\ref{subsec:hd} (see \eqref{treth}),
\ben
\Big[\frac{u-c_m}{u-c_1}
\prod_{a=1}^{m-1}\frac{(u-c_a)^2}{(u-c_a)^2-1}
\Big]_{u=c_m}=\frac{H(\mu)}{H(\la)}.
\een
Furthermore, observe that
\ben
N+\frac{(2\tss c_m-1)(c_m-c_1)}{2\tss c_m}=\frac{(c_m+c_1)(2\tss c_m+1)}{2\tss c_m}.
\een
Hence, to complete the proof of \eqref{trete} we need to verify the identity
\begin{multline}
\frac{(c_m+c_1)(2\tss c_m+1)}{2\tss c_m}
\ts\prod_{a=1}^{m-1}\frac{(c_m+c_a)^2-1}{(c_m+c_a)^2}
\\
{}=\prod_{(i,j)\in\la}\big(N-1+d_{\la}(i,j)\big)\big/
\prod_{(i,j)\in\mu}\big(N-1+d_{\mu}(i,j)\big),
\label{idvery}
\end{multline}
where $d_{\la}(i,j)$ and $d_{\mu}(i,j)$ denote the parameters $d(i,j)$ associated
with the diagrams $\la$ and $\mu$, respectively. The diagram $\la$ is obtained from
$\mu$ by adding one box. Let $(k,l)$ be this box so that $l=\la_k$ and $c_m=\la_k-k+c_1$.
First consider the case $k\leqslant l$. The product on the left hand side
does not depend on the standard tableau $U$ and only depends on its shape $\mu$.
Therefore, the product can be written as
\ben
\prod_{a=1}^{m-1}\frac{(c_m+c_a)^2-1}{(c_m+c_a)^2}
=\prod_{(i,j)\in\mu}\frac{(c_m+c(i,j))^2-1}{(c_m+c(i,j))^2},
\een
where $c(i,j)=j-i+(N-1)/2$ is the content of the box $(i,j)$. We split the
product into two parts by multiplying all terms corresponding to
the subset of boxes $(i,j)$ with $i<l$ and those corresponding to the
subset of boxes $(i,j)$ with $i\geqslant l$. After canceling common factors,
the first part of the product will take the form
\begin{multline}
\frac{c_m+c(l-1,1)-1}{c_m+c(1,1)}\ts\prod_{i=1}^{l-1}
\frac{c_m+c(i,\mu_i)+1}{c_m+c(i,\mu_i)}\\
{}=\frac{c_m+c(l-1,1)-1}{c_m+c(1,1)}\ts\frac{(2\tss c_m)^2}{(2\tss c_m)^2-1}
\ts\prod_{i=1}^{l-1}
\frac{c_m+c(i,\la_i)+1}{c_m+c(i,\la_i)},
\non
\end{multline}
while the second part of the product can be written as
\ben
\frac{c_m+c(l,\la_l)+1}{c_m+c(l,1)}\ts\prod_{j=1}^{\mu_l}
\frac{c_m+c(\la'_j,j)-1}{c_m+c(\la'_j,j)}.
\een
Therefore, the left hand side of \eqref{idvery} equals
\beql{lhsve}
\frac{2\tss c_m\tss (c_m+c(l,\la_l)+1)}{2\tss c_m-1}\ts
\prod_{i=1}^{l-1}
\frac{c_m+c(i,\la_i)+1}{c_m+c(i,\la_i)}\ts
\prod_{j=1}^{\mu_l}
\frac{c_m+c(\la'_j,j)-1}{c_m+c(\la'_j,j)}.
\eeq
To see that this coincides with the right hand side of \eqref{idvery}, note that
for most of the pairs $(i,j)$
the corresponding factors in the numerator
and denominator of the fraction on the right hand side of \eqref{idvery} cancel.
The remaining pairs are divided
into five types: $(i,k)$ with $1\leqslant i<k$; $(k,k)$; $(k,j)$ with
$k<j<l$; $(k,l)$; and $(l,j)$ with $1\leqslant j\leqslant \mu_l$.
Examining the factors for
each of the five types we conclude that their product coincides with \eqref{lhsve}.
Indeed, taking the first type of pairs as an illustration, we obtain
the product
\ben
\prod_{i=1}^{k-1}\ts\frac{N-1+d_{\la}(i,k)}{N-1+d_{\mu}(i,k)}
=\prod_{i=1}^{k-1}\ts\frac{N-1+\la_i+\la_k-i-k+1}{N-1+\la_i+\la_k-i-k}
=\prod_{i=1}^{k-1}\ts
\frac{c_m+c(i,\la_i)+1}{c_m+c(i,\la_i)}
\een
which agrees with the corresponding factors in \eqref{lhsve}. The calculations
for the other types and
the argument in the case $k>l$ are quite similar and will be omitted.
This concludes the proof of the first part of the proposition.
By \eqref{idemptrbd}, the dimension $\dim L(\la)$ equals the trace
$\tr_{1,\dots,m} E_T$ so that the second part
follows from the first by an obvious induction.
\epf

Now we turn to the symplectic group $\Spr_N$, $N=2n$, acting in the space
of tensors \eqref{tenprk} and the commuting action of the Brauer algebra $\Bc_m(-N)$.
The action of $\Bc_m(-N)$ in the space \eqref{tenprk}
is now defined by
\beql{braactsp}
s_{ab}\mapsto -P_{ab},\qquad \ep_{ab}\mapsto -Q_{ab},\qquad
a< b,
\eeq
where $P_{ab}$ is defined in \eqref{pdef}, and
\beql{qdefsp}
Q_{ab}=\sum_{i,j=1}^N \ve_i\tss\ve_j\ts 1^{\ot(a-1)}\ot e_{ij}
\ot 1^{\ot(b-a-1)}\ot e_{i'j'}\ot 1^{\ot(m-b)}.
\eeq

We use Lemma~\ref{lem:idium} in the same way
as for the orthogonal group and write down a recurrence relation for the
respective functions $A_m(u)$ defined by
\eqref{amuf} with the definition \eqref{jmdef} of the Jucys--Murphy elements.
It takes the form
\ben
\bal
\Big(1- \frac{1}{(u-x_{m-1})^2}\Big)A_m(u)
 =
\Big(1- \frac{1}{(u+x_{m-1})^2}\Big)A_{m-1}(u)
-\frac{2\ts(2\tss u-1)\ts x_{m-1}}{(u-x_{m-1})^2(u+x_{m-1})^2}.
\eal
\een
Noting that $A_1(u)=N/(u-c_1)$ with $c_1=(-N-1)/2$
and
using the substitution
\ben
A_m(u)=
\wt{A}_m(u)+ \frac{2u-1}{2u}
\een
we come to the solution
\ben
A_m(u)=\prod_{a=1}^{m-1}\frac{(u-x_a)^2}{(u-x_a)^2-1}
\prod_{a=1}^{m-1}\frac{(u+x_a)^2-1}{(u+x_a)^2}
\Big(\frac{2\tss N}{2\tss u+N+1}-\frac{2\tss u-1}{2\tss u}\Big)+\frac{2\tss u-1}{2\tss u}.
\een

For any diagram $\rho$ with at most $n$ rows set
\beql{dimfsp}
D(\rho)=\frac{1}{H(\rho)}\prod_{(i,j)\in\rho}\big(N+1+d(i,j)\big),
\eeq
where the parameters $d(i,j)$ are now defined by
\ben
d(i,j)=\begin{cases} \rho_i+\rho_j-i-j+1&\qquad
\text{if} \quad i> j,\\
-\rho^{\ts\prime}_i-\rho^{\ts\prime}_j+i+j-1&\qquad\text{if} \quad i\leqslant j.
\end{cases}
\een
The following proposition recovers the symplectic version
of the hook dimension formula \cite[(3.29)]{ek:di}. We suppose that $T$ is a standard
tableau of shape $\la\vdash m$ such that the first row of $\la$ does not exceed $n$,
and $U$ is the tableau obtained from $T$ by removing the box occupied by $m$.
The diagram $\mu$ is the shape of $U$.

\bpr\label{prop:partsp}
We have the relation
\beql{tretesp}
\tr_m E_T=E_U\ts \frac{D(\la')}{D(\mu')}.
\eeq
Moreover,
the dimension of the irreducible representation $L(\rho)$ of $\Spr_N$
equals $D(\rho)$.
\epr

\bpf
Applying again \eqref{jmform}, we get
\ben
\tr_m E_T=E_U\Big[(u-c_m)A_m(u)\Big]_{u=c_m}
\een
so that by the above formula for $A_m(u)$ we have
\ben
\tr_m E_T=\frac{(2\tss c_m+1)(-c_m-c_1)}{2\tss c_m}\ts
E_U\Big[\frac{u-c_m}{u-c_1}
\prod_{a=1}^{m-1}\frac{(u-c_a)^2}{(u-c_a)^2-1}
\Big]_{u=c_m}\prod_{a=1}^{m-1}\frac{(c_m+c_a)^2-1}{(c_m+c_a)^2}.
\een
As with the orthogonal case (Proposition~\ref{prop:parto}),
the proof is reduced to verifying the identity
\begin{multline}
\frac{(-c_m-c_1)(2\tss c_m+1)}{2\tss c_m}
\ts\prod_{a=1}^{m-1}\frac{(c_m+c_a)^2-1}{(c_m+c_a)^2}
\\
{}=\prod_{(i,j)\in\la}\big(N+1+d_{\la'}(i,j)\big)\big/
\prod_{(i,j)\in\mu}\big(N+1+d_{\mu'}(i,j)\big),
\non
\end{multline}
where $d_{\la'}(i,j)$ and $d_{\mu'}(i,j)$ denote the respective parameters $d(i,j)$ associated
with the diagrams $\rho=\la'$ and $\rho=\mu'$.
This identity holds because after the replacement of $N$ by $-N$
it turns into \eqref{idvery}, while the latter can be regarded
as an identity of rational functions
in a variable $N$.
Finally, the dimension of $L(\la')$ equals
$\tr_{1,\dots,m} E_T$ by \eqref{idemptrc}, so that
an obvious induction yields $\dim L(\la')=D(\la')$.
\epf

\section{Images of central idempotents}
\label{sec:ici}
\setcounter{equation}{0}

We consider the cases of orthogonal group $\Or_N$ and
symplectic group $\Spr_N$ (the latter with $N=2n$), simultaneously, unless stated otherwise.
Suppose that $\la$ is a diagram with $m$ boxes such that
$\la'_1\leqslant n$ in the orthogonal case and $\la_1\leqslant n$
in the symplectic case. Define the respective normalized central idempotents $\phi_{\la}$
by \eqref{philao} and \eqref{philasp} and regard them as elements of the algebra
\eqref{tenpendrk} under the action of the Brauer algebra defined by
\eqref{braact} and \eqref{braactsp}. We aim to calculate the images $\ch(\phi_{\la})$
of $\phi_{\la}$ under the characteristic maps \eqref{charabcd}.

\subsection{Main theorem}
\label{subsec:mt}

We let $Y$ run over the Cartan subalgebra of the Lie algebra $\g_N$
so that $Y$ is a diagonal matrix
\ben
Y=\diag(y_1,\dots,y_n,-y_n,\dots,-y_1)\OR
Y=\diag(y_1,\dots,y_n,0,-y_n,\dots,-y_1)
\een
for $N=2n$ or $N=2n+1$, respectively.

Consider the map $F:\g_N\to \Gr_N$
\cite[Theorem~5.2]{oo:ss2}, defined
in a neighborhood of $0$ by the formula
\ben
F(Y)=1+Y^2/2+Y(1+Y^2/4)^{1/2}.
\een
We let $t$ be a complex variable and let $Z=Z(t)$ be the image of the
matrix $tY$ under this map,
\ben
Z=\diag(z_1,\dots,z_n,z^{-1}_n,\dots,z^{-1}_1)\OR
Z=\diag(z_1,\dots,z_n,1,z^{-1}_n,\dots,z^{-1}_1),
\een
respectively.
In particular, writing $Z=F(tY)$ as a power series in $t$ we have
the following first few terms
\ben
Z=1+t\tss Y+\frac12\ts t^2\tss Y^2+\dots.
\een
Therefore, $\ch(\phi_{\la})$ will be found from the coefficient of $t^m$
in the power series expansion
\beql{philach}
\tr\ts \phi_{\la}\ts (Z_1-1)\dots (Z_m-1)=m!\ts t^m\ts \ch(\phi_{\la})+\dots,
\eeq
where the trace is taken over all $m$ copies of $\End\CC^N$ in
\eqref{tenpendrk}. We have
\ben
\tr\ts \phi_{\la}\ts (Z_1-1)\dots (Z_m-1)
=\sum_{k=0}^m (-1)^{m-k}\ts \sum_{a_1<\dots<a_k}
\tr\ts\phi_{\la}\ts Z_{a_1}\dots Z_{a_k}.
\een
Each product $Z_{a_1}\dots Z_{a_k}$ can be written as $P\ts Z_1\dots Z_k\ts P^{-1}$,
where $P$ is the image in \eqref{tenpendrk} of a permutation $p\in\Sym_m$
such that $p(r)=a_r$ for $r=1,\dots,k$. Since $\phi_{\la}$ is
proportional to a central idempotent,
it commutes with $P$, and
by the cyclic property of trace we can bring the above expression to the form
\ben
\tr\ts \phi_{\la}\ts (Z_1-1)\dots (Z_m-1)
=\sum_{k=0}^m (-1)^{m-k}\binom{m}{k} \ts
\tr\ts\phi_{\la}\ts Z_1\dots Z_k.
\een
Propositions~\ref{prop:parto} and \ref{prop:partsp} imply the formula for
the partial trace,
\ben
\tr_m\ts \phi_{\la}=\sum_{\mu}\phi_{\mu},
\een
summed over the diagrams $\mu$ obtained from $\la$ by removing one box.
Hence for any value of the parameter
$k=0,\dots,m$ we have the formula for multiple partial traces
taken over the copies $k+1,\dots,m$ of $\End\CC^N$,
\ben
\tr^{}_{k+1,\dots,m}\ts\phi_{\la}=\sum_{\mu\vdash k,\ts \mu\subseteq\la}\ts
\dim\la/\mu\ts\ts \phi_{\mu},
\een
where, as before, $\dim\la/\mu$ is the number of standard
tableaux with entries in $\{k+1,\dots,m\}$ of the skew shape
$\la/\mu$. Therefore,
\beql{trzmone}
\tr\ts \phi_{\la}\ts (Z_1-1)\dots (Z_m-1)
=\sum_{k=0}^m\ \sum_{\mu\vdash k,\ts \mu\subseteq\la}
(-1)^{|\la|-|\mu|}\dim\la/\mu\ts\binom{|\la|}{|\mu|} \ts
\tr^{}_{1,\dots,k}\ts\phi_{\mu}\ts Z_1\dots Z_k.
\eeq

On the other hand, by \eqref{idemptrbd}
and \eqref{idemptrc},
\ben
\tr^{}_{1,\dots,k}\ts\phi_{\mu}\ts Z_1\dots Z_k
=\frac{\dim\mu}{D(\mu)}\ts
\chi^{\oa_N}_{\mu}(z_1,\dots,z_n)
\een
in the orthogonal case, and
\ben
\tr^{}_{1,\dots,k}\ts\phi_{\mu}\ts Z_1\dots Z_k
=\frac{\dim\mu}{D(\mu')}\ts
\chi^{\spa_N}_{\mu'}(z_1,\dots,z_n)
\een
in the symplectic case. Now, using the notation
\eqref{cnu}, \eqref{doubsf} and \eqref{amun}, we apply
the binomial formula
of \cite[Theorem~1.2]{oo:ss2} which gives
\beql{chirho}
\frac{\chi_{\rho}(z_1,\dots,z_n)}{D(\rho)}
=\sum_{\nu}\frac{s_{\nu}(a_{\rho}\ts |\ts a)\ts
s_{\nu}(t^2y_1^2,\dots,t^2y_n^2)}{C(\nu)},
\eeq
summed over partitions $\nu$ of length not exceeding $n$,
where $\chi_{\rho}(z_1,\dots,z_n)$ denotes any one of the characters
$\chi^{\oa_N}_{\rho}(z_1,\dots,z_n)$ or $\chi^{\spa_N}_{\rho}(z_1,\dots,z_n)$.

This formula implies that if $m$ odd, then the coefficient of $t^m$
on the right hand side of \eqref{trzmone} is zero. Now we assume that
$m$ is even, $m=2l$. Then the coefficient of $t^{2l}$ in the right hand side of \eqref{chirho}
can only come from the terms with the partition $\nu$ having exactly $l$ boxes.
Hence, using \eqref{philach} and \eqref{trzmone} we find that
$\ch(\phi_{\la})$ is the linear combination of the Schur polynomials
$s_{\nu}(y_1^2,\dots,y_n^2)$ with $\nu\vdash l$ occurring with the respective coefficients
\ben
\frac{1}{C(\nu)}
\sum_{\mu\subseteq\la}(-1)^{|\mu|}\ts\ts
\frac{\dim\mu}{|\mu|!}\ts\frac{\dim\la/\mu}{(|\la|-|\mu|)!}\ts s_{\nu}(a_{\rho}\ts |\ts a),
\een
where $\rho=\mu$ and $\rho=\mu'$ in the orthogonal and symplectic case, respectively.
Thus, recalling the notation $H(\theta)$ in \eqref{prhook} we arrive at the main result.

\bth\label{thm:mainch}
Suppose that $\la$ is a diagram with $m$ boxes such that
$\la'_1\leqslant n$ in the orthogonal case and $\la_1\leqslant n$
in the symplectic case. Then the image $\ch(\phi_{\la})$ of the normalized central
idempotent $\phi_{\la}$ under the respective characteristic map \eqref{charabcd}
is zero if $m$ is odd. If $m=2l$ is even, then the image is found by
\beql{chlao}
\ch(\phi_{\la})
=\sum_{\nu\vdash l}\frac{s_{\nu}(y_1^2,\dots,y_n^2)}{C(\nu)}
\sum_{\mu\subseteq\la}(-1)^{|\mu|}\ts
\frac{s_{\nu}(a_{\mu}|\ts a)}{H(\mu)\ts H(\la/\mu)}
\eeq
in the orthogonal case, and
\beql{chlasp}
\ch(\phi_{\la})
=\sum_{\nu\vdash l}\frac{s_{\nu}(y_1^2,\dots,y_n^2)}{C(\nu)}
\sum_{\mu\subseteq\la}(-1)^{|\mu|}\ts
\frac{s_{\nu}(a_{\mu'}|\ts a)}{H(\mu)\ts H(\la/\mu)}
\eeq
in the symplectic case.
\qed
\eth

Note that by the vanishing theorem \cite{o:qi} we have $s_{\nu}(a_{\rho}\ts|\ts a)=0$
unless $\nu\subseteq\rho$. Therefore, the first sum in \eqref{chlao}
is restricted to the partitions $\nu$ contained in $\la$, while the second sum
is restricted to the partitions $\mu$ containing $\nu$. Similarly,
the first sum in \eqref{chlasp}
is restricted to the partitions $\nu$ contained in $\la'$, while the second sum
is restricted to the partitions $\mu$ containing $\nu'$.

\bex\label{ex:twobytwo}
Consider the orthogonal case with $\la=(2^2)$. By \eqref{chlao},
the image $\ch(\phi_{(2^2)})$ is a linear combination for
the Schur polynomials $s_{\nu}(y_1^2,\dots,y_n^2)$ with $\nu=(2)$ and $\nu=(1^2)$.
Using \cite[Proposition~3.2]{ms:lr}, we find
\ben
\bal
s_{(2)}(a_{(2)}\ts|\ts a)&=(a_{n+2}-a_n)(a_{n+2}-a_{n+1}),\\
s_{(2)}(a_{(2\ts 1)}\ts|\ts a)&=(a_{n+2}-a_{n-1})(a_{n+2}-a_{n+1}),\\
s_{(2)}(a_{(2^2)}\ts|\ts a)&=(a_{n+2}-a_{n-1})(a_{n+2}-a_{n})
\eal
\een
and
\ben
\bal
s_{(1^2)}(a_{(1^2)}\ts|\ts a)&=(a_{n+1}-a_{n-1})(a_{n}-a_{n-1}),\\
s_{(1^2)}(a_{(2\ts 1)}\ts|\ts a)&=(a_{n+2}-a_{n-1})(a_{n}-a_{n-1}),\\
s_{(1^2)}(a_{(2^2)}\ts|\ts a)&=(a_{n+2}-a_{n-1})(a_{n+1}-a_{n-1}).
\eal
\een
For the sequence $a_i=(\ve+i-1)^2$ we have $a_{n+i}-a_{n+j}=(i-j)(N+i+j-2)$.
Hence the sums in \eqref{chlao} are found by
\ben
\frac{s_{(2)}(a_{(2)}\ts|\ts a)}{4}-\frac{s_{(2)}(a_{(2\ts 1)}\ts|\ts a)}{3}
+\frac{s_{(2)}(a_{(2^2)}\ts|\ts a)}{12}=1
\een
and
\ben
\frac{s_{(1^2)}(a_{(1^2)}\ts|\ts a)}{4}-\frac{s_{(1^2)}(a_{(2\ts 1)}\ts|\ts a)}{3}
+\frac{s_{(1^2)}(a_{(2^2)}\ts|\ts a)}{12}=1.
\een
Thus,
\ben
\bal
\ch(\phi_{(2^2)})&=\frac{1}{(N-1)N(N+1)(N+2)}\ts s_{(2)}(y_1^2,\dots,y_n^2)\\
{}&+\frac{1}{(N-3)(N-2)(N-1)N}\ts s_{(1^2)}(y_1^2,\dots,y_n^2).
\eal
\een
\eex

\subsection{Symmetrizers and antisymmetrizers}
\label{subsec:sa}

Now we consider the particular cases, where $\la$ is a row or column diagram
with $2l$ boxes. In each of these cases there is a unique standard tableau $T$
of shape $\la$ so that by \eqref{philao} and \eqref{philasp}, $\phi_{\la}$
is proportional to the primitive
idempotent $E_T$. For $\la=(2l)$ the primitive idempotent coincides with the
{\it symmetrizer\/} $S^{(2l)}$, while for $\la=(1^{2l})$ it coincides with
the {\it antisymmetrizer\/} $A^{(2l)}$ in the
Brauer algebra. We will produce
the images $\ch(S^{(2l)})$ and $\ch(A^{(2l)})$
in an explicit form.
Suppose first that $\la=(2l)$ in the orthogonal case. Then the first sum in
\eqref{chlao} contains only one term with $\nu=(l)$, while the
second sum is taken over row-diagrams
$\mu=(k)$ with $l\leqslant k\leqslant 2l$. By \eqref{doubsf}
we have
\ben
\bal
s_{(l)}(x\ts|\ts a)&=\sum_{i_1\leqslant\dots\leqslant i_l}
(x_{i_1}-a_{i_1})\dots (x_{i_l}-a_{i_l+l-1})\\
{}&=\sum_{i_1\leqslant\dots\leqslant i_l}
(x_{n-i_1+1}-a_{i_1})\dots (x_{n-i_l+1}-a_{i_l+l-1}),
\eal
\een
where the second relation holds since $s_{(l)}(x\ts|\ts a)$
is a symmetric polynomial. Recalling the definition \eqref{amun},
we find $a_{(k)}=(a_{k+n},a_{n-1},\dots,a_1)$. Hence, taking $x=a_{(k)}$
we find that the only nonzero summand corresponds to $i_1=\dots=i_l=n$,
\beql{slrow}
s_{(l)}(a_{(k)}\ts|\ts a)=(a_{k+n}-a_n)\dots (a_{k+n}-a_{n+l-1}).
\eeq
Furthermore, recalling that $a_i=(\ve+i-1)^2$ we find
\ben
\bal
s_{(l)}(a_{(k)}\ts|\ts a)
&=\big((k+n+\ve-1)^2-(n+\ve-1)^2\big)\\
{}&\times\dots\times \big((k+n+\ve-1)^2-(n+l+\ve-2)^2\big)
=\frac{k!\ts (N+k+l-3)!}{(k-l)!\ts (N+k-3)!}.
\eal
\een
The sum in \eqref{chlao} then equals
\ben
\bal
&\sum_{k=l}^{2l}(-1)^k\ts
\frac{(N+k+l-3)!}{(k-l)!\ts (N+k-3)!\ts(2l-k)!}\\
{}&=\sum_{r=0}^l (-1)^{l-r}\binom{N+2l-3}{l-r}\ts
\binom{N+2l-3+r}{r}=1.
\eal
\een
Thus, taking into account the constants $D(\la)$ for $\la=(2l)$ and $C(\nu)$
for $\nu=(l)$ we come
to the following corollary (for a different proof see
\cite[Proposition~3.4]{m:ff}).

\bco\label{cor:rowo}
The image of the symmetrizer $S^{(2l)}\in\Bc_{2l}(N)$ under the characteristic
map is found by
\ben
\ch(S^{(2l)})=\frac{N+4l-2}{(2l)!\ts(N+2l-2)}\ts
\sum_{1\leqslant i_1\leqslant\dots\leqslant i_l\leqslant n}
y_{i_1}^2\dots y_{i_l}^2.
\een
\eco

Now let
$\la=(1^{2l})$ with $2l\leqslant n$.
The second sum in
\eqref{chlao} is now taken over column-diagrams
$\mu=(1^k)$ with $l\leqslant k\leqslant 2l$ and $\nu=(1^l)$.
Using \cite[Proposition~3.2]{ms:lr}, we find that
\beql{scova}
s_{(1^l)}(a_{(1^k)}\ts|\ts a)=(a_{n-l+2}-a_{n-k+1})\dots(a_{n+1}-a_{n-k+1}).
\eeq
Under the specialization $a_i=(\ve+i-1)^2$ this simplifies to
\ben
s_{(1^l)}(a_{(1^k)}\ts|\ts a)=\frac{k!\ts (N-k)!}{(k-l)!\ts (N-k-l)!}
\een
so that the sum in \eqref{chlao} equals
\ben
\sum_{k=l}^{2l}(-1)^{k}\ts
\frac{(N-k)!}{(k-l)!\ts (2l-k)!\ts(N-k-l)!}=(-1)^l,
\een
thus leading to the image of the antisymmetrizer.

\bco\label{cor:colo}
The image of the antisymmetrizer $A^{(2l)}\in\Bc_{2l}(N)$ under the characteristic
map is found by
\ben
\ch(A^{(2l)})=\frac{(-1)^l}{(2l)!}\ts
\sum_{1\leqslant i_1<\dots< i_l\leqslant n}
y_{i_1}^2\dots y_{i_l}^2.
\een
\eco
Note that this result also follows easily from the observation that
$A^{(2l)}$ coincides with the antisymmetrizer in the group algebra
$\CC[\Sym_{2l}]$. Indeed, it suffices to apply \eqref{idemptr} with $\la=(1^{2l})$
and replace $X$ by the diagonal matrix $Y$.

The calculation in the symplectic case is quite similar.
Suppose first that $\la=(2l)$ with $2l\leqslant n$.
Then $\nu=(1^l)$ in \eqref{chlasp} and $\mu$ runs over diagrams $(k)$ with
$l\leqslant k\leqslant 2l$. Using \eqref{scova} with the sequence $a_i=i^2$
and performing the same calculations as in the orthogonal case
we find the image of $S^{(2l)}$; see also \cite[Proposition~3.5]{m:ff}.

\bco\label{cor:rowsp}
The image of the symmetrizer $S^{(2l)}\in\Bc_{2l}(-2n)$ under the characteristic
map is found by
\ben
\ch(S^{(2l)})=\frac{(-1)^l\ts (n-2l+1)}{(2l)!\ts(n-l+1)}\ts
\sum_{1\leqslant i_1<\dots< i_l\leqslant n}
y_{i_1}^2\dots y_{i_l}^2.
\een
\eco

Finally, if $\la=(1^{2l})$ then $\nu=(l)$ in \eqref{chlasp}
and $\mu$ runs over diagrams $(1^k)$ with
$l\leqslant k\leqslant 2l$. Applying now \eqref{slrow},
we calculate the image of $A^{(2l)}$.

\bco\label{cor:colsp}
The image of the antisymmetrizer $A^{(2l)}\in\Bc_{2l}(-2n)$ under the characteristic
map is found by
\ben
\ch(A^{(2l)})=\frac{1}{(2l)!}\ts
\sum_{1\leqslant i_1\leqslant\dots\leqslant i_l\leqslant n}
y_{i_1}^2\dots y_{i_l}^2.
\een
\eco
This result also follows from the observation that
$A^{(2l)}$ coincides with the symmetrizer in the group algebra
$\CC[\Sym_{2l}]$. It suffices to apply \eqref{idemptr} with $\la=(2l)$
and replace $X$ by the diagonal matrix $Y$.

\end{document}